# TRAINING SAMPLES IN OBJECTIVE BAYESIAN MODEL SELECTION[1]

BY JAMES O. BERGER AND LUIS R. PERICCHI[2]

*Duke University and Universidad de Puerto Rico*

Central to several objective approaches to Bayesian model selection is the use of training samples (subsets of the data), so as to allow utilization of improper objective priors. The most common prescription for choosing training samples is to choose them to be as small as possible, subject to yielding proper posteriors; these are called minimal training samples.

When data can vary widely in terms of either information content or impact on the improper priors, use of minimal training samples can be inadequate. Important examples include certain cases of discrete data, the presence of censored observations, and certain situations involving linear models and explanatory variables. Such situations require more sophisticated methods of choosing training samples. A variety of such methods are developed in this paper, and successfully applied in challenging situations.

**1. Introduction.** Training samples play a central role in a variety of statistical methodologies, including classification and discrimination, cross-validation, robustness and model selection, from both Bayesian and frequentist perspectives. Two recent developments in Bayesian model selection are the *intrinsic Bayes factor* of Berger and Pericchi (1996a) and the *expected posterior prior* of Pérez (1998) and Pérez and Berger (2002). Central to both is utilization of training samples to convert improper objective priors into the proper distributions typically needed for model selection. The most common prescription for choosing training samples is to choose them to be

Received July 2002; revised April 2003.

[1]Supported by NSF Grants DMS-98-02261 and DMS-01-03265, and by grants from the School of Natural Sciences, UPR-RRP and CONICIT-USB-Venezuela.

[2]Supported by a Guggenheim Fellowship.

*AMS 2000 subject classifications.* Primary 62F03, 62F15; secondary 62N03, 62B10, 62F40.

*Key words and phrases.* Intrinsic Bayes factors, expected posterior priors, training samples, objective priors, intrinsic priors, censored data, linear models.

as small as possible, subject to yielding proper posteriors; these are called *minimal training samples*.

While fine for many problems, minimal training samples have been found to be suboptimal in an ever-increasing number of important statistical situations, in particular those in which the data can vary widely in terms of information content. Important examples include the presence of censored observations, studied in Section 3; certain cases of discrete data, studied in Section 4; and situations involving unbalanced linear models or covariates, studied in Section 5.

A variety of strategies have been developed to overcome the limitation of minimal training samples, and the main purpose of this paper is to outline these strategies. The generalizations of training samples considered herein can alternatively be viewed as choosing training samples in a random fashion, or as providing a "weighting" to chosen training samples. One particularly interesting example is a *sequential random minimal training sample*, which is a training sample of smallest size such that the posterior is proper, but which is obtained by drawing observations randomly, without replacement, from the set of data. Another natural use of random training samples is when the original data is not available, but sufficient statistics are given; training samples can then be generated from the conditional distribution of the data, given the sufficient statistics.

We will see considerable evidence that use of the new definitions of training samples can successfully overcome a wide variety of problems in Bayesian model selection. It is worth noting up front, however, that we were unable to define any type of "optimal" training sample; the paper can thus be viewed as providing a useful set of strategies that can be employed to obtain good training samples, with statistical judgement being required to select from among these strategies in particular contexts. While this prevents the proposed model selection methods from being completely automatic, the judgements involved in choosing good training samples will typically be much less than the judgements needed to implement an actual subjective Bayesian analysis. See Section 6 for overall suggestions and further context concerning this issue.

In the remainder of this section, the model selection problem is stated, and intrinsic Bayes factors and expected posterior priors are defined. Section 1.3 discusses the key problem that arises, which can be best understood through the device of studying the *intrinsic priors* corresponding to intrinsic Bayes factors; these are the priors that, if used directly to compute Bayes factors, would yield (in an asymptotic sense) the same answers as the intrinsic Bayes factors. As further discussed in Berger and Pericchi (2001), we feel this to be a powerful unifying approach to understanding the performance of default Bayes factors.

There has been a significant literature discussing training samples in these and other Bayesian contexts. Other recent articles include Gelfand, Dey and Chang (1992), de Vos (1993), Iwaki (1997, 1999), Lingham and Sivaganesan (1997, 1999), Alqallaf and Gustafson (2001) and Ghosh and Samanta (2002).

1.1. *Model selection notation.* Suppose that we are comparing $q$ models for the data $\mathbf{X} = (X_1, \ldots, X_n)$,

$$M_i : \mathbf{x} \text{ has density } f_i(\mathbf{x}|\boldsymbol{\theta}_i), \qquad i = 1, \ldots, q,$$

where the $\boldsymbol{\theta}_i$ are unknown model parameters. Let $\pi_i(\boldsymbol{\theta}_i)$, $i = 1, \ldots, q$, be prior distributions for the unknown parameters, and define the marginal or predictive densities of $\mathbf{x}$,

$$m_i(\mathbf{x}) = \int f_i(\mathbf{x}|\boldsymbol{\theta}_i) \pi_i(\boldsymbol{\theta}_i) \, d\boldsymbol{\theta}_i.$$

The *Bayes factor* of $M_j$ to $M_i$ is given by

$$(1) \qquad B_{ji} = \frac{m_j(\mathbf{x})}{m_i(\mathbf{x})} = \frac{\int f_j(\mathbf{x}|\boldsymbol{\theta}_j) \pi_j(\boldsymbol{\theta}_j) \, d\boldsymbol{\theta}_j}{\int f_i(\mathbf{x}|\boldsymbol{\theta}_i) \pi_i(\boldsymbol{\theta}_i) \, d\boldsymbol{\theta}_i}$$

and is often interpreted as the "odds provided by the data for $M_j$ versus $M_i$." Thus $B_{ji} = 10$ would suggest that the data favor $M_j$ over $M_i$ at odds of ten to one. Alternatively, $B_{ji}$ is sometimes called the "weighted likelihood ratio of $M_j$ to $M_i$," with the priors being the "weighting functions." These interpretations are particularly appropriate when, as here, we focus on conventional or default choices of the priors.

1.2. *Intrinsic Bayes factors and expected posterior priors.* For the $q$ models $M_1, \ldots, M_q$ suppose that only noninformative priors $\pi_i^{\mathrm{N}}(\boldsymbol{\theta}_i)$, $i = 1, \ldots, q$, are available. In general, we recommend that these be chosen to be "reference priors" [see Berger and Bernardo (1992)]. Define the corresponding marginal or predictive densities of $\mathbf{x}$,

$$m_i^{\mathrm{N}}(\mathbf{x}) = \int f_i(\mathbf{x}|\boldsymbol{\theta}_i) \pi_i^{\mathrm{N}}(\boldsymbol{\theta}_i) \, d\boldsymbol{\theta}_i.$$

Unfortunately, the direct use of improper priors for defining Bayes factors in (1) is not generally justifiable [cf. Berger and Pericchi (1996a, 2001)], but they can be utilized for model selection through the introduction of training samples. Here is the standard type of training sample.

DEFINITION 0 [Berger and Pericchi (1996a)]. A training sample, to be indexed by $l$, is a subset of the data, $\mathbf{x}(l)$. It is called *proper* if $0 < m_i^{\mathrm{N}}(\mathbf{x}(l)) < \infty$ for all $M_i$. Let $\mathcal{X}^{\mathrm{P}}$ denote the set of all proper training samples and define its cardinality as $L_{\mathrm{P}}$. A training sample is *minimal* if it is proper and no subset is proper. A *minimal training sample* will be denoted MTS; let $\mathcal{X}^{\mathrm{M}}$ and $L_{\mathrm{M}}$ denote, respectively, the set of all MTS and its cardinality.

Thus $\mathbf{x}(l)$ can be used to "convert" the improper $\pi_i^{\mathrm{N}}(\boldsymbol{\theta}_i)$ to proper posteriors,

$$\pi_i^{\mathrm{N}}(\boldsymbol{\theta}_i|\mathbf{x}(l)) = \frac{f_i(\mathbf{x}(l)|\boldsymbol{\theta}_i)\pi_i^{\mathrm{N}}(\boldsymbol{\theta}_i)}{m_i^{\mathrm{N}}(\mathbf{x}(l))}. \tag{2}$$

These posteriors can then be used to define Bayes factors for the remaining data.

Since there are typically many possible training samples, it is natural to average the resulting Bayes factors over the training samples in some fashion. The resulting Bayes factor for comparing $M_j$ to $M_i$ [called the *intrinsic Bayes factor* (IBF) in Berger and Pericchi (1996a)] is

$$B_{ji} = B_{ji}^{\mathrm{N}} \cdot \mathrm{AVE}[B_{ij}^{\mathrm{N}}(\mathbf{x}(l))], \tag{3}$$

where

$$B_{ji}^{\mathrm{N}} = B_{ji}^{\mathrm{N}}(\mathbf{x}) = \frac{m_j^{\mathrm{N}}(\mathbf{x})}{m_i^{\mathrm{N}}(\mathbf{x})} \quad \text{and} \quad B_{ij}^{\mathrm{N}}(l) = B_{ij}^{\mathrm{N}}(\mathbf{x}(l)) = \frac{m_i^{\mathrm{N}}(\mathbf{x}(l))}{m_j^{\mathrm{N}}(\mathbf{x}(l))},$$

and "AVE" denotes an average of the $B_{ij}^{\mathrm{N}}(\mathbf{x}(l))$. A variety of possible averages have been considered [see Berger and Pericchi (1996a, 2001)], the most common being arithmetic, geometric and median averages. Some recent references to use and development of intrinsic Bayes factors in various scenarios include Berger and Pericchi (1996b, 1996c, 1998), Bertolino and Racugno (1996), De Santis and Spezzaferri (1997), Lingham and Sivaganesan (1997, 1999), Sun and Kim (1997), Berger, Pericchi and Varshavsky (1998), Key, Pericchi and Smith (1999), Moreno, Bertolino and Racugno (1998, 1999, 2001), Bertolino, Racugno and Moreno (2000), Berger and Mortera (1999), Sivaganesan and Lingham (1999), Kim and Sun (2000), Rodriguez and Pericchi (2001), Beattie, Fong and Lin (2002), Ghosh and Samanta (2002) and Paulo (2002).

Another recent use of training samples for model selection is in the development of *empirical expected posterior priors* [Pérez (1998), Pérez and Berger (2001, 2002) and Neal (2001)], defined as

$$\pi_i^{\mathrm{EP}}(\boldsymbol{\theta}_i) = \frac{1}{L_{\mathrm{M}}} \sum_{\mathbf{x}(l) \in \mathcal{X}^{\mathrm{M}}} \pi_i^{\mathrm{N}}(\boldsymbol{\theta}_i|\mathbf{x}(l)). \tag{4}$$

The idea is that, instead of using the minimal training samples to define proper posteriors for computation of Bayes factors and then averaging the ensuing Bayes factors, one can first average the proper posteriors and then compute Bayes factors with the results. This approach can be embedded within Markov chain Monte Carlo analysis, which can be a considerable computational advantage. Another advantage is that one can use minimal training samples for each separate model, which has certain computational and theoretical benefits.

1.3. *Evaluation of intrinsic Bayes factors and a key condition.* The most basic approach to evaluation of intrinsic Bayes factors is simply to see if they produce sensible answers. In Berger and Pericchi (1996c, 2001) it is argued that the best way to study this is to determine the *intrinsic prior* corresponding to an IBF. The intrinsic prior is that prior which would yield Bayes factors that are approximately equal to the IBF, in an asymptotic sense. If this intrinsic prior is sensible, then the IBF is judged to be sensible. The power and sensitivity of the use of intrinsic priors in appraising default Bayesian model selection methods is illustrated in Berger and Mortera (1999) and Berger and Pericchi (2001); see also the Examples in Sections 3 and 4 in this paper. It is particularly important to establish the existence (and sensibility) of intrinsic priors when new concepts are introduced (as here, to deal with censored data and other difficulties); such initial study can give considerable confidence that the new IBFs will work more generally.

One can also use intrinsic priors directly as the conventional prior for model selection [cf. Sun and Kim (1997), Moreno, Bertolino and Racugno (1998, 1999, 2001), Bertolino, Racugno and Moreno (2000) Kim and Sun (2000), Cano, Kessler and Moreno (2002), Moreno, Girón and Torres (2004), Moreno, Torres and Casella (2002), Paulo (2002), Girón, Martínez and Moreno (2003) and Moreno and Liseo (2003)]; this is an attractive possibility, although it is often more computationally intensive than using the IBF directly. Indeed, analytic determination of intrinsic priors can itself be quite difficult, and they will frequently not have closed form expressions. [They can have expressions amenable to MCMC computation, however; see Pérez and Berger (2002).]

Computation of intrinsic priors corresponding to model selection requires an extension from the finite set of proper training samples for the existing data to a hypothetical sampling space of proper training samples, to be denoted by $\mathcal{X}^{\mathrm{I}}$, based on imagining availability of an infinite sequence of data. Choice of this sampling space is sometimes automatic, but sometimes involves judgement; an example of each is given below. Note that $\mathcal{X}^{\mathrm{I}}$ will typically be considered fixed for all models under consideration, although there are situations (such as with expected posterior priors) in which $\mathcal{X}^{\mathrm{I}}$ can be allowed to vary with the model.

EXAMPLE 1. Suppose $X_1, X_2, \ldots$ are i.i.d. from the normal distribution with unknown mean $\mu$ and variance $\sigma^2$. For the usual reference prior, $\pi(\mu, \sigma^2) = 1/\sigma^2$, an easy computation shows that an MTS must consist of any two distinct observations. Thus, if we use the MTS notion to define training samples, it is clear that we should define $\mathcal{X}^{\mathrm{I}}$ to be the set of all pairs of (distinct) observations from the hypothetical infinite population of normal observations having mean $\mu$ and variance $\sigma^2$. (The word "distinct" is theoretically superfluous, since the distribution is absolutely continuous.)

EXAMPLE 2. Consider a linear model in which observation $x_i$ has associated $k$-vector of covariates $\mathbf{D}_i$, $i = 1, \ldots, n$. Suppose that an MTS would consist of any $m$ observations for which the corresponding vectors $\mathbf{D}_i$ are linearly independent. If we wish to extend this definition to an infinite population, it is necessary to decide if the covariates are viewed as fixed or themselves random. In the former case, we can simply imagine that the hypothetical infinite population arises from proportionally replicated covariates. Letting $\mathbf{D}$ denote the $n \times k$ design matrix of fixed covariates, $\mathcal{X}^{\mathrm{I}}$ can then be formally defined as the space of sets of $m$ observations that arise by first randomly drawing $m$ linearly independent rows from $\mathbf{D}$, and then generating corresponding observations from the linear model. If the covariates are considered random, one would first have to define the sampling distribution of covariates and then construct $\mathcal{X}^{\mathrm{I}}$ by draws from the covariate distribution, followed by generation of observations from the linear model. In this paper we shall only consider the fixed covariates scenario.

The special case of intrinsic priors that will be considered in this paper is that in which there are two models, $M_0$ nested in $M_1$, and the arithmetic average is used in (3). Then the intrinsic prior is given by

(5) $$\pi_1^{\mathrm{I}}(\boldsymbol{\theta}_1) = \pi_1^{\mathrm{N}}(\boldsymbol{\theta}_1) E_{\boldsymbol{\theta}_1}^{M_1}[B_{01}^{\mathrm{N}}(\mathbf{X}(l)) | \mathcal{X}^{\mathrm{I}}],$$

where $E^{M_1}$ refers to expectation under model $M_1$. This expression differs from the earlier expressions for an intrinsic prior, given in Berger and Pericchi (1996a, 1996c), because of the conditioning on $\mathcal{X}^{\mathrm{I}}$. In the examples considered in these earlier papers, $P_{\theta_1}^{M_1}(\mathcal{X}^{\mathrm{I}}) = 1$, so that the conditioning was not needed. In general, however, the conditioning is needed to correctly define the intrinsic prior.

One important property of a "good" intrinsic prior is that it integrate to one. [If it fails to do so, the corresponding IBF would appear to be "biased" toward one of the models; see, e.g., Berger and Mortera (1999) and Berger and Pericchi (2001).] Theorem 1 in Berger and Pericchi (1996a) asserts that this will be so (under mild regularity conditions) if $\pi_0^{\mathrm{N}}$ is proper (trivially satisfied if $M_0$ is a simple hypothesis). Again, however, it was implicitly assumed that $\mathcal{X}^{\mathrm{I}}$ had probability one; in this paper we formally state our assumption:

ASSUMPTION 0.　$P_{\boldsymbol{\theta}_i}^{M_i}(\mathcal{X}^{\mathrm{I}}) = 1$, $i = 0, 1$.

In Sections 3 and 4, we will see that this assumption can be violated for the set of minimal training samples, in situations involving censoring or when inappropriate initial noninformative priors are utilized.

If Assumption 0 is satisfied and $\pi_0^{\mathrm{N}}$ is proper, then the intrinsic prior will be proper. For simplicity we only show this in the case when $M_0$ is a simple model.

LEMMA 1. *If Assumption* 0 *holds,* $M_0$ *is a simple model* (*i.e.,* $\boldsymbol{\theta}_0$ *is specified*) *and the intrinsic prior is given by* (5), *then*

$$\int \pi_1^{\mathrm{I}}(\boldsymbol{\theta}_1) \, d\boldsymbol{\theta}_1 = 1.$$

PROOF. Since, by Assumption 0, $\mathcal{X}^{\mathrm{I}}$ is the support of $\mathbf{X}(l)$ under $M_1$, and since $\mathcal{X}^{\mathrm{I}}$ contains only proper training samples, it follows from (5) that

$$\int \pi_1^{\mathrm{I}}(\boldsymbol{\theta}_1) \, d\boldsymbol{\theta}_1 = \int\int_{\mathcal{X}^{\mathrm{I}}} \pi_1^{\mathrm{N}}(\boldsymbol{\theta}_1) \frac{m_0^{\mathrm{N}}(\mathbf{x}(l))}{m_1^{\mathrm{N}}(\mathbf{x}(l))} f_1(\mathbf{x}(l)|\boldsymbol{\theta}_1) \, d\mathbf{x}(l) \, d\boldsymbol{\theta}_1.$$

Applying Fubini's theorem to switch the order of integration yields

$$\int \pi_1^{\mathrm{I}}(\boldsymbol{\theta}_1) \, d\boldsymbol{\theta}_1 = \int_{\mathcal{X}^{\mathrm{I}}} m_0^{\mathrm{N}}(\mathbf{x}(l)) \, d\mathbf{x}(l) = P_{\boldsymbol{\theta}_0}^{M_0}(\mathcal{X}^{\mathrm{I}}) = 1,$$

the last step following from the assumption that $M_0$ is simple and Assumption 0. □

If Assumption 0 does not hold, the intrinsic prior can be highly unsatisfactory (even improper, as we will see in later examples), casting considerable doubt on the quality of the associated IBF. Thus, if Assumption 0 is violated in a particular context, the set of training samples should be enlarged until the assumption is satisfied. This can sometimes be done by changing the noninformative prior but, more generally, a more sophisticated definition of training sample is required.

Note that under Assumption 0 the intrinsic prior in (5) has the alternative representation

(6) $$\pi_1^{\mathrm{I}}(\boldsymbol{\theta}_1) = \int_{\mathcal{X}^{\mathrm{I}}} \pi_1^{\mathrm{N}}(\boldsymbol{\theta}_1|\mathbf{x}(l)) m_0^{\mathrm{N}}(\mathbf{x}(l)) \, d\mathbf{x}(l),$$

which is also called the *base-model posterior expected prior* in Pérez and Berger (2002). If one is interested in utilizing the intrinsic prior directly in computing Bayes factors, this expression is typically most useful in that, within MCMC, one can simply drop the integral sign and treat $\mathbf{x}(l)$ as a latent variable. The improved training samples that are obtained in the following sections for IBFs can also be immediately utilized in (6) to obtain improved intrinsic priors that are computationally attractive.

As a final comment, when $\pi_0^{\mathrm{N}}(\boldsymbol{\theta}_0)$ is improper, then $\pi_1^{\mathrm{I}}(\boldsymbol{\theta}_1)$ will also be improper. However, it is well calibrated with $\pi_0^{\mathrm{N}}(\boldsymbol{\theta}_0)$, in the sense that a limiting argument over compact sets shows that the Bayes factor for the two priors is a well-defined limit of proper priors. See Berger and Pericchi (1996a) for discussion and Moreno, Bertolino and Racugno (1998) for implementation.

**2. Generalizations of training samples.** To handle situations in which Assumption 0 is violated and in which training samples can contain very different information, it is necessary to introduce more general types of training samples.

2.1. *Randomized and weighted training samples.*

DEFINITION 1. A *randomized training sample* with sampling mechanism $\mathbf{p} = (p_1, \ldots, p_{L_P})$, where $\mathbf{p}$ is a probability vector, is obtained by drawing a training sample from $\mathcal{X}^P$ according to $\mathbf{p}$. Alternatively, the training samples can be considered to be *weighted training samples* with weights $p_i$.

EXAMPLE 3 (Sequential random sampling). We will be particularly interested in sequential minimal training samples (SMTS) that are each obtained by drawing observations from the collection of data $\mathbf{x} = \{x_1, x_2, \ldots, x_n\}$ by simple random sampling (without replacement for a given SMTS), stopping when the subset so formed, $\mathbf{x}^*(l) = (x(l)_1, \ldots, x(l)_{N(l)})$, is a proper training sample. Note that $N(l)$ is itself a random variable. Although intuitively and operationally one obtains an SMTS by sequential random sampling, such training samples can also be described via Definition 1, with $p_i$ being the probability of obtaining the $i$th SMTS via sampling without replacement from the set of observations, and all other proper training samples being assigned probability 0.

REMARK. When the $X_i$ are i.i.d. and arise from an absolutely continuous distribution, then an SMTS will typically equal an MTS with probability one, since each distinct observation will typically have the same effect on posterior propriety.

EXAMPLE 4 (Sampling of minimal training samples). Often the number of minimal training samples $L_M$ is extremely large, so that the computation of the averages in (3) can be very expensive. In such situations it usually suffices to just randomly choose minimal training samples [i.e., set $p_i = 1/L_M$ for $\mathbf{x}(l) \in \mathcal{X}^M$ and set $p_i = 0$ otherwise in Definition 1]. Indeed, in Varshavsky (1995) the theory of $U$-statistics is used to indicate that it often suffices to randomly choose $L = kn$ minimal training samples, where $n$ is the sample size of the actual data and $k$ is the size of the minimal training sample (assuming there is a fixed size). This is clearly much smaller than the number of minimal training samples, $\binom{n}{k}$. (Unfortunately, precise guidelines as to the choice of $L$ are not available, so a reasonable practical implementation is to start with the choice $kn$ and increase $L$ until the change in the resulting Bayes factor is sufficiently small.)

EXAMPLE 5 (Probability proportional to information). Observations are often associated with covariates. In linear models for instance, a training sample $\mathbf{x}(l)$ will typically have a corresponding "design matrix" of covariates $\mathbf{D}(l)$ and corresponding "information" proportional to $|\mathbf{D}(l)'\mathbf{D}(l)|$. One could choose training samples with probability proportional to this information (or perhaps the square root of the information). This was proposed in de Vos (1993).

On the other hand, one does not want training samples to be too informative. Suppose, for instance, that almost all of the information in the entire sample is due to a single observation. Utilization of that observation as a training sample can be inappropriate, as will be seen in Section 5. Indeed, it is generally a good idea to restrict attention to training samples that contain only a modest fraction of the total information in the data, although this may not always be possible [cf. Rodriguez and Pericchi (2001)].

EXAMPLE 6 (Random sampling to reach a given information level). An interesting variant of the sequential random sampling approach to construction of a training sample is to stop, not when the training sample is proper, but when the training sample contains a certain amount of "information." We do not pursue this idea here.

2.2. *Imaginary training samples.* A different notion that has been employed [in, e.g., Good (1950), Smith and Spiegelhalter (1980), Iwaki (1997, 1999), Pérez (1998), Rodriguez and Pericchi (2001), Ghosh and Samanta (2002) and Pérez and Berger (2002)] is that of an imaginary training sample: training samples are generated, not from the real data, but from some specified distribution. For instance, in model selection one might elicit a subjective predictive distribution, $m^*(\mathbf{x}^*)$, where $\mathbf{x}^*$ is thought of as a "future" minimal training sample. One could then draw training samples from this distribution for Bayesian model selection, or use the associated expected posterior priors [see Pérez and Berger (2002) for motivation and further discussion].

One potential difficulty with training sample methods is that often only sufficient statistics (and not the actual data) are available. Use of imaginary training samples can overcome this difficulty.

DEFINITION 2. A *conditional imaginary training sample*, for a situation in which only sufficient statistics from a model are available, is defined to be a training sample from the conditional distribution of the data given the sufficient statistics.

If $S$ is a sufficient statistic, the factorization theorem gives
$$f(\mathbf{x}|\theta) = g(S|\theta) \cdot h(\mathbf{x}|S),$$

and we can repeatedly draw conditional imaginary training samples $\mathbf{x}^*$ from the corresponding marginal distribution $h(\mathbf{x}^*|S)$. In computation of intrinsic Bayes factors or expected posterior priors one then presumes that the imaginary training sample $\mathbf{x}^*$ arose from the density $f(\mathbf{x}^*|\theta)$.

EXAMPLE 7 (Example 1 continued). Let $X_1,\ldots,X_n$ be an i.i.d. sample from the normal distribution with mean $\mu$ and variance $\sigma^2$, but suppose that only the sufficient statistics $\bar{x}$ and $s^2 = \sum_i (x_i - \bar{x})^2$ are reported, along with $n$. A very simple way to draw conditional imaginary training samples is to create the surrogate data set

$$X_i^* = (Z_i - \bar{Z})\frac{s}{s_Z} + \bar{x}, \qquad i = 1,\ldots,n,$$

where the $Z_i$ are independent standard normal with sample mean and sum of squared deviations $\bar{Z}$ and $s_Z^2$, respectively. This surrogate data set clearly has the same sample mean and sum of squared deviations as the original data and is a draw from $h(\mathbf{x}|\bar{x}, s^2)$. One can then choose training samples (recall minimal training samples were of size 2) from this surrogate data set. (Note that it is necessary to have $n \geq 3$ in order to have training samples that are not simply the entire data set.) One can also draw additional surrogate data sets if more training samples are needed (an advantage of using imaginary training samples). Imaginary training samples are used as if they were real training samples, that is, they are assumed to arise from the original normal distribution with $\mu$ and $\sigma^2$.

EXAMPLE 8 (Poisson distribution). Suppose that $X$ is a single realization from a Poisson distribution with mean $\theta T$, arising as the number of rare events observed in a time period $T$. We consider testing of $H_0 : \theta = \theta_0$ versus $H_1 : \theta \neq \theta_0$, utilizing the improper Jeffreys prior, $\pi_1^N(\theta) = \theta^{-1/2}$.

A natural way to define imaginary training samples is to use the fact that such a Poisson $X$ can be viewed as arising from a sum of the indicators of events occurring with exponential inter-arrival times. More precisely, for $i = 1,\ldots$, consider $X_i \sim f(x_i|\theta) = \theta\exp(-\theta x_i)$, and define

$$X \equiv \left\{\text{first } j \text{ such that } S_j = \sum_{i=1}^j X_i > T\right\} - 1.$$

Then $X$ has the Poisson distribution with mean $\theta T$.

It is natural to utilize these latent $\{x_1,\ldots,x_X\}$ to construct imaginary training samples. No simple trick is available as in the previous example, so we must determine $h(x_1,\ldots,x_X|X)$. Computation yields that this is the uniform density on $\sum_{i=1}^X x_i < T$. Thus, if training samples consist of a single observation (as is the case in the testing situation we consider with the Jeffreys prior), an imaginary training sample can be drawn from the marginal distribution of a single $x_i$ arising from this uniform distribution, which is

$$h(x_i|X) = \frac{X}{T}\left(1 - \frac{x_i}{T}\right)^{X-1}, \qquad 0 < x_i < T. \tag{7}$$

Single imaginary training samples can thus be drawn as $X_i^* = T[1 - U^{1/X}]$, where $U$ is Uniform$(0,1)$. These are then used in constructing intrinsic Bayes factors and/or expected posterior priors, as if they had arisen from the exponential density with mean $1/\theta$. Note that we have implicitly assumed that $T > 0$ in defining the imaginary training samples.

The situation is not always as nice as the above examples would suggest, in that the information needed to construct $h(\mathbf{x}|S)$ in order to generate the imaginary training samples can be lost when a sufficiency reduction is effected.

EXAMPLE 9 (Linear model). Suppose $\mathbf{Y}(n \times 1)$ arises from the linear model

$$\mathbf{Y} = \mathbf{X}\boldsymbol{\beta} + \varepsilon, \qquad \varepsilon \sim \mathcal{N}_n(\mathbf{0}, \sigma^2 \mathbf{I}_n),$$

where $\boldsymbol{\beta} = (\beta_1, \beta_2, \ldots, \beta_k)'$ is unknown, $\sigma^2$ is known, and $\mathbf{X}$ is an $(n \times k)$ given design matrix of rank $k \leq n$. The least squares estimate $\hat{\boldsymbol{\beta}} = (\mathbf{X}'\mathbf{X})^{-1}\mathbf{X}'\mathbf{y}$ is then sufficient for $\boldsymbol{\beta}$, and one might be presented only with $n$, $\hat{\boldsymbol{\beta}}$ and its covariance matrix $\boldsymbol{\Sigma} = \sigma^2(\mathbf{X}'\mathbf{X})^{-1}$ after a sufficiency reduction. From this one cannot reconstruct the conditional distribution of the data given $\hat{\boldsymbol{\beta}}$, because the design matrix, and hence the covariates, have been "lost" (unless $n = k$, in which case it can be reconstructed from $\boldsymbol{\Sigma}$). So imaginary training samples cannot be generated in this way. For some ideas as to alternative ways of generating imaginary training samples in situations such as this, see Iwaki (1999).

2.3. *Utilization of generalized training samples.* For training samples defined as in Definition 1 and considered as weighted training samples the arithmetic IBF and empirical expected posterior priors are defined, respectively, as

$$B_{ji}^{\text{A}} = B_{ji}^{\text{N}} \sum_{l=1}^{L_{\text{P}}} p_l B_{ij}^{\text{N}}(\mathbf{x}(l)), \tag{8}$$

$$\pi_i^{\text{EP}}(\boldsymbol{\theta}_i) = \sum_{l=1}^{L_{\text{P}}} p_l \pi_i^{\text{N}}(\boldsymbol{\theta}_i|\mathbf{x}(l)). \tag{9}$$

It is often not feasible to compute these weighted averages (because of the large number of possible training samples), in which case it is easier to draw $L$ random training samples, $\mathbf{x}(1), \mathbf{x}(2), \ldots, \mathbf{x}(L)$, according to the random schemes discussed above for generating the training samples (repeats allowed), and then just approximate the arithmetic IBF and empirical expected posterior priors by, respectively,

$$B_{ji}^{\mathrm{A}} \cong B_{ji}^{\mathrm{N}} \frac{1}{L} \sum_{l=1}^{L} B_{ij}^{\mathrm{N}}(\mathbf{x}(l)), \tag{10}$$

$$\pi_i^{\mathrm{EP}}(\boldsymbol{\theta}_i) \cong \frac{1}{L} \sum_{l=1}^{L} \pi_i^{\mathrm{N}}(\boldsymbol{\theta}_i | \mathbf{x}(l)). \tag{11}$$

It usually suffices to take $L$ to be a modest multiple of the overall sample size $n$.

**3. Censoring.** Censored data provides a key illustration of these ideas. We begin with an example involving right-censoring. For another discussion of training samples in the presence of censoring, see Lingham and Sivaganesan (1999).

EXAMPLE 10 (Right censoring of exponential data). Suppose the data $x_1, \ldots, x_n$ arises as a random sample from the right-censored Exponential$(\theta)$ density; thus, if $x_i < r$ it arises from the density $f(x_i|\theta) = \theta \exp(-\theta x_i)$, while $P(X_i = r|\theta) \equiv p(\theta) = \exp(-r\theta)$. It is desired to test

$$M_0 : \theta = \theta_0 \quad \text{versus} \quad M_1 : \theta \neq \theta_0.$$

Consider the usual default prior $\pi_1^{\mathrm{N}}(\theta) = \theta^{-1}$. It is easy to show that any single uncensored observation yields a proper posterior, while no number of censored observations will do so. Hence the set of minimal training samples $\mathcal{X}^{\mathrm{M}}$ consists of the collection of single uncensored observations. Since censored observations never enter into the training samples, the MTS's will intuitively be biased in favor of larger values of $\theta = 1/E(X_i|\theta)$, which seems undesirable.

To evaluate the situation more carefully, consider the intrinsic prior for $\theta$ corresponding to the arithmetic IBF; this prior is given by (5), where the sampling space of training samples, here denoted $\mathcal{X}^{\mathrm{MI}}$, is simply the interval $(0, r)$ [i.e., the space of single uncensored observations drawn from $f(x|\theta, x < r)$]. Note first that Assumption 0 is violated, since

$$P_{\theta_i}^{M_i}(\mathcal{X}^{\mathrm{MI}}) = P_{\theta_i}^{M_i}(X < r) = 1 - \exp(-r\theta_i) < 1, \qquad i = 0, 1,$$

so that we expect problems with the intrinsic prior (and hence with the intrinsic Bayes factor). Noting that the Bayes factor for a training sample is $B_{01}^{N}(x) = \theta_0 \exp(-\theta_0 x) / \int \frac{1}{\theta} \theta \exp(-\theta x) \, d\theta = x\theta_0 \exp(-\theta_0 x)$, the intrinsic prior in (5) is given by

$$\pi^{I}(\theta) = \frac{1}{\theta} \int_0^r x\theta_0 \exp(-\theta_0 x) \frac{\theta \exp(-\theta x)}{(1 - \exp(-r\theta))} \, dx$$

$$= \frac{\theta_0}{(1 - \exp(-r\theta))} \left[ \frac{1}{(\theta + \theta_0)^2} - e^{-(\theta + \theta_0)r} \left( \frac{r}{\theta + \theta_0} + \frac{1}{(\theta + \theta_0)^2} \right) \right].$$

This is not a proper prior; indeed, as $\theta \to 0$ the prior behaves like a constant times $1/\theta$, which is nonintegrable, a particularly egregious failing.

One possible solution to this problem would be to use a noninformative prior that enlarges the set of MTS's. Indeed, for this problem involving right censoring the Jeffreys-rule prior is $\pi^{J}(\theta) = \theta^{-1}[1 - \exp(-r\theta)]^{1/2}$ [De Santis, Mortera and Nardi (2001)]. For this prior, it can be shown that any single observation, censored or uncensored, is an MTS, so that Assumption 0 is trivially satisfied and the resulting intrinsic prior must integrate to one. Note, however, that extra work is involved in finding the Jeffreys-rule prior, and this can be formidable in more complex situations (e.g., in Example 11). Furthermore, the intrinsic prior that results from use of the Jeffreys-rule prior here has the quite unappealing property (see the Appendix) that its median is $O(r^{-1})$ as $r \to 0$. This unattractive behavior arises because the highly informative training samples (the uncensored observations) have effects averaged with the (many more) censored observations that have negligible information content as $r \to 0$. Hence we turn to use of sequential minimal training samples to solve the problem.

For the prior $\pi^{N}(\theta) = \theta^{-1}$ a SMTS is of the form $\mathbf{x}(l) = (r, \ldots, r, x(l))$, where $x(l)$ is the first uncensored observation that arises in simple random sampling (without replacement) from the data. (In contrast, none of the $r$ would be present in an MTS.) The natural sampling space for such training samples is the set $\mathcal{X}^{SI}$ of possible sequences $\mathbf{x}(l) = (r, \ldots, r, x(l))$ of i.i.d observations arising from the censored exponential distribution. Let $N^*(l) = N(l) - 1$ denote the number of censored observations in the SMTS from $\mathcal{X}^{SI}$, and write $p(\theta) = P(X > r|\theta) = \exp(-\theta r)$. Note that $P(N^*(l) = j|\theta) = (1 - p(\theta))p(\theta)^j$, and that the joint density of $\mathbf{x}(l)$ is $f(\mathbf{x}(l)|\theta) = p(\theta)^j \theta \exp(-\theta x(l))$.

Letting $n_u$ denote the number of uncensored observations in the actual data, and letting $T$ denote the sum of all observations (censored and uncensored), computation yields

$$B_{10}^{N} = \frac{m_1^{N}(\mathbf{x})}{m_0^{N}(\mathbf{x})} = \Gamma(n_u)(T\theta_0)^{-n_u} e^{T\theta_0},$$

$$B_{01}^{\mathrm{N}}(\mathbf{x}(l)) = \frac{m_0^{\mathrm{N}}(\mathbf{x}(l))}{m_1^{\mathrm{N}}(\mathbf{x}(l))} = \theta_0(N^*(l)r + x(l))\exp(-[N^*(l)r + x(l)]\theta_0).$$

The approximate arithmetic IBF in (10), corresponding to $L$ random SMTS draws, is then

$$B_{10}^{\mathrm{A}} \cong \Gamma(n_u)(T\theta_0)^{-n_u} e^{T\theta_0} \frac{1}{L} \sum_{l=1}^{L} \theta_0(N^*(l)r + x(l))\exp(-[N^*(l)r + x(l)]\theta_0).$$

To investigate the behavior of this IBF, we again study its corresponding intrinsic prior. From (5) and noting that $P_{\theta_i}^{M_i}(\mathcal{X}^{\mathrm{SI}}) = 1$, this is given by

$$
\begin{aligned}
\pi^{\mathrm{I}}(\theta) &= \frac{1}{\theta} E_\theta^{M_1}[B_{01}^{\mathrm{N}}(\mathbf{x}^*(l))] \\
(12) \quad &= \frac{1}{\theta} \sum_{j=0}^{\infty} \int_0^r \theta_0(jr + x)\exp(-[jr + x]\theta_0)p(\theta)^j \theta \exp(-\theta x)\,dx \\
&= \frac{\theta_0}{(\theta + \theta_0)^2},
\end{aligned}
$$

the last step following from standard calculations involving geometric series. This is a very sensible intrinsic prior for the problem, being proper and having median equal to $\theta_0$. Indeed, this is the intrinsic prior for the exponential testing problem when no censoring is present and ordinary MTS are used [Pericchi, Fiteni and Presa (1993)], an appealing result. The indication is that use of SMTS leads to a very satisfactory arithmetic IBF in the presence of censoring.

It would be fascinating if the result observed in Example 10—that the intrinsic prior in the presence of censoring and using SMTS equals the intrinsic prior when there is no censoring and using MTS—held in general. Unfortunately, this is not the case, as can be seen by considering the density $f(x|\theta) = (0.5)\exp(-|x - \theta|)$, together with a constant default prior on $\theta$. Detailed calculations yield that the intrinsic prior without censoring and using MTS is not equal to the intrinsic prior with right censoring and RMTS. We omit the details.

That the intrinsic prior in (12) is proper would not have needed exact calculation. Indeed, consider the general case of censoring of i.i.d. observations, with a *known* censoring mechanism and the use of SMTS. Then the natural sampling space is the set $\mathcal{X}^{\mathrm{SI}}$ of possible sequences of i.i.d observations arising from the original distribution (with censoring), with the sampling stopping the first time the training sample is proper. Assuming that the sampling is guaranteed to stop with probability one for any of the models and parameter values under consideration (i.e., that the sampling mechanism is a *proper stopping rule*), then Assumption 0 is satisfied and Lemma 1 shows that the intrinsic prior is proper.

When the censoring mechanism is at least partly unknown, intrinsic priors cannot be defined. However, SMTS can be defined, and the corresponding IBFs or empirical expected posterior priors utilized to compute Bayes factors. We illustrate this with an example comparing two exponential distributions.

EXAMPLE 11 (Comparison of two exponential populations). The following data, which appeared in Gehan (1965), were analyzed in Cox and Oakes (1984) as arising from (possibly censored) exponential distributions. The data show times of remission (as measured by freedom from symptoms), in weeks, of leukemia patients, where the first group consists of control individuals and the second group consists of individuals treated with the drug 6-mercaptopurine. The data is as follows, where + indicates that the data has been censored.

*Control*: 1, 1, 2, 2, 3, 4, 4, 5, 5, 8, 8, 8, 8, 11, 11, 12, 12, 15, 17, 22, 23.
*Treated*: 6+, 6, 6, 6, 7, 9+, 10+, 10, 11+, 13, 16, 17+, 19+, 20+, 22, 23, 25+, 32+, 32+, 34+, 35+.

Notice that the control group has no censored observations, but more than half of the observations from the treated set have been censored.

Following Cox and Oakes (1984) and with $j = 1, 2$ referring to the control and treatment groups, respectively, assume that the uncensored failure times $t_{ji}$ follow the Exponential($\theta_j$) distribution. Write each observation as $x_{ji} = (y_{ji}, v_{ji})$, where $y_{ji} = \min(t_{ji}, c_{ji})$, with $c_{ji}$ denoting the censoring time (known for the actual data), and $v_{ji} = 0$ if $t_{ji} \leq c_{ji}$ (uncensored) and $v_{ji} = 1$ otherwise. Specifying the density here is problematical when the overall distribution of the $c_{ji}$ is not known, but for Bayesian analysis we only need the likelihood function of $(\theta_1, \theta_2)$ for the given data, and this is given by

$$(13) \qquad \text{lik}(\theta_1, \theta_2) = \prod_{j=1}^{2} \left[ \prod_{i=1}^{n_{ju}} \theta_j e^{-t_{ji}\theta_j} \prod_{i=1}^{n_{jc}} e^{-t_{ji}\theta_j} \right],$$

where $n_{ju}$ and $n_{jc}$ denote, respectively, the number of uncensored and censored observations in each group, and the labels are rearranged if necessary.

We want to test the hypotheses

$$(14) \qquad M_0 : \theta_1 = \theta_2 = \theta \quad \text{versus} \quad M_1 : \theta_1 \neq \theta_2.$$

In the analysis, we will utilize the usual noninformative priors $\pi_0^N(\theta) = \theta^{-1}$ and $\pi_1^N(\theta_1, \theta_2) = \theta_1^{-1} \theta_2^{-1}$. As in Example 10, it then follows that an SMTS must consist of a sequence of censored observations from each group, followed by an uncensored observation. (Since in the actual data the control group contains only uncensored observations, an SMTS for this data will contain just a single control observation, but we will write down expressions for the general case.) Write this SMTS as

$$\mathbf{x}(l) = \begin{cases} c_{11}(l), \ldots, c_{1N_1^*}(l), t_{11}(l), \\ c_{21}(l), \ldots, c_{2N_2^*}(l), t_{21}(l), \end{cases}$$

where $N_1^*$ and $N_2^*$ are the (random) stopping times in obtaining the SMTS.

Straightforward calculation then shows the arithmetic IBF to be

$$(15) \qquad B_{10}^A = \frac{\Gamma(n_{1u})\Gamma(n_{2u})}{\Gamma(n_{1u} + n_{2u})} \frac{(T_1 + T_2)^{n_{1u}+n_{2u}}}{T_1^{n_{1u}} T_2^{n_{2u}}} \frac{1}{L} \sum_{l=1}^{L} \frac{T_1(l) T_2(l)}{(T_1(l) + T_2(l))^2},$$

where

$$T_1 = \sum_{i=1}^{n_{1u}} t_{1i} + \sum_{j=1}^{n_{1c}} c_{1j}, \qquad T_2 = \sum_{i=1}^{n_{2u}} t_{2j} + \sum_{j=1}^{n_{2c}} c_{2j},$$

$$T_1(l) = \sum_{j=1}^{N_1^*(l)} c_{1j}(l) + t_{11}(l), \qquad T_2(l) = \sum_{j=1}^{N_2^*(l)} c_{2j}(l) + t_{21}(l)$$

and $L$ is the number of SMTS that are to be drawn.

For analysis of the actual data above we computed (15) using $L = n = 42$, $L = 2n$ and $L = 5n$ training samples obtained by simple random sampling (without replacement) from the data. The resulting Bayes factors were $B_{10} = 544$, 493 and 584, respectively, showing decisive evidence against the null model and only modest variation with respect to the number of training samples drawn. If equal prior probabilities are assumed for the hypotheses, then the posterior probability of $M_1$ is about $P(M_1|\mathbf{x}) = 0.998$.

It is also straightforward to calculate the approximations to the empirical expected posterior priors, given in (11), and use them to compute the Bayes factor of $M_1$ to $M_0$. The result is

$$B_{10}^{\text{EP}} = \frac{\Gamma(n_{1u}+1)\Gamma(n_{2u}+1)}{\Gamma(n_{1u}+n_{2u}+2)}$$

$$\times \frac{\sum_{l=1}^{L} T_1(l) T_2(l) (T_1 + T_1(l))^{-(n_{1u}+1)} (T_2 + T_2(l))^{-(n_{2u}+1)}}{\sum_{l=1}^{L} (T_1(l) + T_2(l))^2 (T_1 + T_2 + T_1(l) + T_2(l))^{-(n_{1u}+n_{2u}+2)}}.$$

For the data above and random training samples of sizes $L = n = 42$, $L = 2n$ and $L = 5n$, the resulting Bayes factors were $B_{10}^{\text{EP}} = 742, 713$ and 728, respectively. These are similar to the arithmetic IBF, but are systematically somewhat larger, providing support for the suggestion in Pérez and Berger (2002) that the empirical expected posterior priors will yield Bayes factors that are somewhat more favorable to the more complex model than IBFs or intrinsic priors.

Perhaps the most interesting feature of the above example is that Bayes factors and posterior probabilities could be computed rather easily, without needing to know the nature of the censoring mechanism. In contrast, classical answers typically depend on the (often unknown) censoring mechanism. This is thus an important situation in which the objective Bayesian approach requires significantly less knowledge than a frequentist approach.

Lack of knowledge of the censoring mechanism does preclude computation of the intrinsic prior corresponding to the arithmetic IBF in censoring situations, however; without such knowledge, it is not clear how to define the sampling space for the SMTS, needed for computation of the intrinsic prior. Of course, one might reasonably "cheat" in this situation, using the suggestion from Example 10 that the intrinsic prior for SMTS and in the presence of censoring might well be close to the intrinsic prior for the problem when there is no censoring (and MTS are used). One could then directly use these "approximate" intrinsic priors to compute the Bayes factor.

EXAMPLE 12 (Example 11 continued). An MTS in the uncensored version of this bi-exponential problem would consist of one observation from each of the control and treatment groups. Denoting this MTS by simply $(t_1, t_2)$, the corresponding intrinsic priors are easily seen to be $\pi_0^I(\theta) = \pi_0^N(\theta) = \theta^{-1}$ and

$$\pi_1^I(\theta_1, \theta_2) = \int_0^\infty \int_0^\infty \frac{t_1 t_2}{(t_1 + t_2)^2} \exp(-t_1 \theta_1) \exp(-t_2 \theta_2) \, dt_1 \, dt_2.$$

Combining these intrinsic priors with the likelihood (13) and interchanging order of integration results in the Bayes factor

$$B_{10}^I = \frac{\Gamma(n_{1u} + 1)\Gamma(n_{2u} + 1)}{\Gamma(n_{1u} + n_{2u})} (T_1 + T_2)^{n_{1u} + n_{2u}}$$

(16)
$$\times \int_0^\infty \int_0^\infty \frac{t_1 t_2}{(t_1 + t_2)^2} \frac{1}{(T_1 + t_1)^{n_{1u}+1}(T_2 + t_2)^{n_{2u}+1}} \, dt_1 \, dt_2.$$

For the data of Example 11 numerical computation yields $B_{10}^I = 503$, a value quite close to those obtained with the approximate arithmetic IBF and using SMTS training samples.

Another advantage of having (approximate) intrinsic priors, as above, is that they can be utilized to develop conditional frequentist tests. Indeed, the intrinsic prior above has been utilized in Paulo (2002) to develop optimal conditional frequentist tests for the bi-exponential testing problem.

**4. Discrete examples.** Difficulties with training sample approaches for discrete data have been highlighted in several papers [e.g., Bertolino and Racugno (1996), O'Hagan (1997) and Berger and Pericchi (1998)]. We first revisit one of the more vexing examples, to see if randomized training samples fix the problem.

EXAMPLE 13 (Bernoulli testing). Based on $n$ Bernoulli trials, with $P(X_i = 1|\theta) = \theta = 1 - P(X_i = 0|\theta)$, it is desired to test

$$M_0 : \theta = \theta_0 \quad \text{versus} \quad M_1 : \theta \neq \theta_0.$$

Suppose the improper Haldane prior $\pi_1^N(\theta) = \theta^{-1}(1-\theta)^{-1}$ is utilized to construct an IBF. This is a quite inferior noninformative prior, but it is interesting to see if IBFs can be made robust to poor choices of the initial noninformative prior. Note that for the Haldane prior

$$B_{10}^N = \frac{\Gamma(S)\Gamma(n-S)}{\Gamma(n)\theta_0^S(1-\theta_0)^{(n-S)}},$$

where $S$ is the number of ones in the data.

With the Haldane prior an MTS must consist of precisely one 1 and one 0. (One and only one of each is needed for the resulting posterior to be proper.) Since $P_\theta^{M_i}(\mathcal{X}^M) = 2\theta(1-\theta) < 1$, Assumption 0 is clearly violated and the resulting IBF is again suspect. Indeed, noting that $B_{01}^N(\{0,1\}) = \theta_0(1-\theta_0)$, it is immediate from (5) that the implied intrinsic prior is

(17) $$\pi^I(\theta) = \frac{\theta_0(1-\theta_0)}{\theta(1-\theta)}.$$

This is itself improper—indeed it is simply a constant multiple of the original Haldane prior—and strongly suggests that the IBF for the Haldane noninformative prior and the usual definition of a minimal training sample do not correspond to a sensible Bayes procedure.

An extreme case of this example arises when $\theta_0 = 0$ and the data consists of one 1 and the rest 0. O'Hagan (1997) noted that then $M_0 : \theta_0 = 0$ is wrong with certainty (one cannot observe a 1 under $M_0$), yet the intrinsic Bayes factor will then equal $1/(n-1)$, for $n \geq 2$. The basic problem, in this case, is that $P_{\theta_0}^{M_0}(\mathcal{X}^M) = 0$, an extreme violation of Assumption 0. A single extra 1 ($S = 2$) would solve the problem, making $B_{10} = \infty$ (as it should be), but the behavior of the IBF is indeed disturbing when $S = 1$.

This extreme example is a good test of the effectiveness of SMTS. An SMTS will either be of the form $\mathbf{x}^*(l) = (0, 0, \ldots, 0, 1)$ or $\mathbf{x}^*(l) = (1, 1, \ldots, 1, 0)$; these can obviously be summarized by specifying $N_0$ (the number of zeroes) and $N_1$ (the number of ones), respectively. Noting that $P_\theta(N_0) = (1-\theta)^{N_0}\theta$ and $P_\theta(N_1) = \theta^{N_1}(1-\theta)$, for $N_0, N_1 = 1, 2, \ldots$, it follows that

$$B_{01}^N(N_0) = \frac{\theta_0(1-\theta_0)^{N_0}}{\int_0^1 \theta(1-\theta)^{N_0} \pi^N(\theta)\, d\theta}$$
$$= N_0 \theta_0 (1-\theta_0)^{N_0},$$
$$B_{01}^N(N_1) = N_1(1-\theta_0)\theta_0^{N_1}.$$

To determine the intrinsic prior corresponding to the arithmetic IBF in (5), we first choose $\mathcal{X}^I$ to be the set of training samples, $N_0$ and $N_1$, arising from an infinite series of Bernoulli$(\theta)$ trials. When $0 < \theta_0 < 1$, it is clear that $P_{\theta_i}^{M_i}(\mathcal{X}^I) = 1$, so that the intrinsic prior is

$$\pi_1^I(\theta) = \pi_1^N(\theta) E_\theta^{M_1}[B_{01}^N(\mathbf{X}(l))]$$
$$= \theta^{-1}(1-\theta)^{-1}\left[\theta_0 \theta \sum_{i=1}^{\infty} i[(1-\theta)(1-\theta_0)]^i + (1-\theta_0)(1-\theta)\sum_{i=1}^{\infty} i[\theta\theta_0]^i\right]$$
$$= \theta_0(1-\theta_0)\left[\frac{1}{(1-(1-\theta)(1-\theta_0))^2} + \frac{1}{(1-\theta\theta_0)^2}\right].$$

It can be verified that $\int_0^1 \pi_1^I(\theta)\, d\theta = 1$, so the intrinsic prior is proper. (With slightly less work, this also follows from Lemma 1.) Also, the intrinsic prior is admirably balanced, in the sense that the median is very close to $\theta_0$. [Numerical computation shows that $0.48 < P(\theta < \theta_0) < 0.52$ for all $\theta_0$.] Thus all indications are that the use of the SMTS has corrected the problem caused by the bad initial noninformative prior.

Of course, we needed the condition $0 < \theta_0 < 1$ for the SMTS to work. For the extreme $\theta_0 = 0$ (or the case $\theta_0 = 1$), Assumption 0 remains violated even for the SMTS; indeed, $P_{\theta_0}^{M_0}(\mathcal{X}^P) = 0$ in the extreme cases, so that no set of proper training samples can work. As an indication of the danger in using training sample approaches when Assumption 0 is violated, consider again the situation considered by O'Hagan (1997). The arithmetic IBF, based on use of SMTS for the given data, can be computed to be $B_{10} = (n^2 - n + 2)/[2n(n-1)]$, which while an improvement over $1/(n-1)$, is still not $\infty$, as it should be. Hence even use of SMTS cannot correct the situation when Assumption 0 is violated.

One might wonder if the the training sample solution fails as, say, $\theta_0 \to 0$. This is awkward to discuss in terms of the arithmetic IBF itself, since the sample size would correspondingly need to grow to $\infty$ before a proper training sample could be obtained. We thus look at direct use of the intrinsic prior (18) to see if it yields a satisfactory Bayes factor. Indeed, the resulting Bayes factor is

$$B_{10}^I = \frac{\int_0^1 \theta^S(1-\theta)^{n-S}[(1-(1-\theta)(1-\theta_0))^{-2} + (1-\theta\theta_0)^{-2}]\, d\theta}{\theta_0^{S-1}(1-\theta_0)^{n-S-1}}.$$

For the problematical case $S = 1$, $n \geq 2$ and $\theta_0 \to 0$,

$$B_{10}^{\text{I}} \to \int_0^1 \theta(1-\theta)^{n-1}\left[\frac{1}{\theta^2} + 1\right] d\theta,$$

which is infinite. Thus, for very small $\theta_0$ and the observation $S = 1$, one would properly conclude that the alternative $M_1$ is true.

Next we revisit the Poisson example from Section 2.2, to see the effectiveness of imaginary training samples when only a sufficient statistic is given.

EXAMPLE 14 (Example 8 continued). Recall we are testing $H_0 : \theta = \theta_0$ versus $H_1 : \theta \neq \theta_0$. For the Jeffreys prior $1/\sqrt{\theta}$ under $H_1$ computation yields that the formal Bayes factor is

$$B_{10}^{\text{N}} = \frac{\Gamma(X + 1/2)}{T^{(X+1/2)} \theta_0^X e^{-T\theta_0}}.$$

Recall that we generate imaginary data $x_i^*$ and assume it to be exponential with mean $1/\theta$. A single such observation is a minimal training sample. The arithmetic IBF in (10) is thus given by

$$B_{01}^{\text{A}} = B_{10}^{\text{N}} \frac{1}{L} \sum_{l=1}^{L} \frac{\theta_0}{\Gamma(3/2)}(x_l^*)^{3/2} \exp(-\theta_0 x_l^*).$$

To study the performance of this objective Bayes factor we again determine the corresponding intrinsic prior. Since the $x_i^*$ were actually generated from (7), the intrinsic prior in (5) is given by

$$(18) \quad \pi^{\text{I}}(\theta) = \frac{1}{\sqrt{\theta}} \lim_{T \to \infty} \int_0^T \frac{\theta_0}{\Gamma(3/2)}(x_l^*)^{3/2} \exp(-\theta_0 x_l^*) \frac{X}{T}\left(1 - \frac{x_l^*}{T}\right)^{X-1} dx_l^*.$$

[The intrinsic prior, as defined in Berger and Pericchi (1996a), is based on letting the sample size go to infinity; for the Poisson problem the analogue of this definition is $T \to \infty$.] Since the integrand in (18) is bounded above by $\theta_0 (x_l^*)^{3/2} \exp(-\theta_0 x_l^*)/\Gamma(3/2)$, which is integrable, we may invoke the dominated convergence theorem to take the limit inside the integral. Furthermore,

$$\lim_{T \to \infty} \frac{X}{T}\left(1 - \frac{x_l^*}{T}\right)^{X-1} = \theta \exp(-\theta x_l^*)$$

almost surely, so that

$$\pi^{\text{I}}(\theta) = \frac{1}{\sqrt{\theta}} \int_0^\infty \frac{\theta_0}{\Gamma(3/2)}(x_l^*)^{3/2} \exp(-\theta_0 x_l^*)\theta \exp(-\theta x_l^*) \, dx_l^*$$

$$= \frac{3\theta_0 \sqrt{\theta}}{2(\theta + \theta_0)^{5/2}}.$$

This is a proper prior, and has median approximately equal to (1.7) $\theta_0$, a quite satisfactory prior. Hence the arithmetic IBF based on imaginary training samples arising from a single Poisson observation seems fine.

**5. Information-based training samples in the linear model.** As mentioned in Section 2.1, it is attractive to consider choosing training samples according to their information content. We begin with a classic example demonstrating the need to do this. [A related example can be found in Iwaki (1997).]

EXAMPLE 15 (Findley's example). Findley (1991) demonstrated the inadequacy of BIC in the following situation. Suppose we observe $X_i = d_i\theta + \varepsilon_i$, for $i = 1, \ldots, n$, and that the $\varepsilon_i$ are i.i.d. $\mathcal{N}(0,1)$. It is desired to test

$$M_0 : \theta = 0 \quad \text{versus} \quad M_1 : \theta \neq 0.$$

The standard noninformative prior is $\pi^N(\theta) = 1$, and the corresponding formal Bayes factor is

$$B_{10}^N = \frac{\sqrt{2\pi}}{\|\mathbf{d}\|} \exp\left(\frac{(\sum_{i=1}^n x_i d_i)^2}{2\|\mathbf{d}\|^2}\right),$$

where

$$\|\mathbf{d}\|^2 = \sum_{i=1}^n d_i^2.$$

A minimal training sample is a single observation $x_i$, and

$$B_{01}^N(x_i) = \frac{|d_i|}{\sqrt{2\pi}} \exp\left(-\frac{x_i^2}{2}\right).$$

It follows that the arithmetic IBF is

$$B_{10}^A = \frac{\sqrt{2\pi}}{\|\mathbf{d}\|} \exp\left(\frac{(\sum_{i=1}^n x_i d_i)^2}{2\|\mathbf{d}\|^2}\right) \frac{1}{n} \sum_{i=1}^n \frac{|d_i|}{\sqrt{2\pi}} \exp\left(-\frac{x_i^2}{2}\right).$$

The interesting special case considered by Findley was $d_i = i^{-1/2}$. Then as $n \to \infty$ it is straightforward to show that $\|\mathbf{d}\|^2 = O(\log n)$,

$$\frac{(\sum_{i=1}^n x_i d_i)^2}{\|\mathbf{d}\|^2} = \theta^2 \log n + 2Z\theta\sqrt{\log n} + O(1)$$

and

$$\frac{1}{n}\sum_{i=1}^n \frac{|d_i|}{\sqrt{2\pi}} \exp\left(-\frac{x_i^2}{2}\right) = O(n^{-1/2}),$$

where $Z$ is a standard normal random variable. It follows that
$$B_{10}^{A} = O([n \log n]^{-1/2}) \exp\left(\tfrac{1}{2}\theta^2 \log n + Z\theta\sqrt{\log n}\,\right)$$
$$= O((\log n)^{-1/2} n^{(\theta^2-1)/2} \exp\left(Z\theta\sqrt{\log n}\,\right)).$$

Under $M_0 : \theta = 0$, it is clear that $B_{10}^{A} = O([n \log n]^{-1/2}) \to 0$ as $n \to \infty$; this is fine, as it indicates that $M_0$ is true. But if $M_1$ is true with $\theta^2 < 1$, then $B_{10}^{A} \to 0$ also, which means the arithmetic IBF is then inconsistent, a severe inadequacy. (If $\theta^2 = 1$, the arithmetic IBF is consistent or not, depending on the sign of $Z$, that is, it will be consistent half the time.) The source of the problem (as with the associated inconsistency of BIC, as shown by Findley) is that the observations $x_i$ contain drastically decreasing information, $d_i^2 = i^{-1}$, as $i$ increases.

This thus provides a good test for the idea of weighting the training samples by the amount of information they contain, that is, setting $p_i = d_i^2 / \|\mathbf{d}\|^2$ in Definition 1, and using the corresponding weighted IBF in (8). The resulting Bayes factor, using similar arguments to above, satisfies
$$B_{10}^{A} = \frac{\sqrt{2\pi}}{\|\mathbf{d}\|} \exp\left(\frac{(\sum_{i=1}^{n} x_i d_i)^2}{2\|\mathbf{d}\|^2}\right) \sum_{i=1}^{n} \frac{|d_i|^3}{\|\mathbf{d}\|^2 \sqrt{2\pi}} \exp\left(-\frac{x_i^2}{2}\right)$$
$$= O((\log n)^{-3/2} n^{\theta^2/2} \exp\left(Z\theta\sqrt{\log n}\,\right)).$$

This still goes to 0 under $M_0$ (as it should), but now goes to $\infty$ under $M_1$ (as it should).

The use of weighted training samples solved the inconsistency problem, but that is a very crude criterion and the goal in use of training samples is to achieve actual Bayesian behavior. Unfortunately, even the use of weighted training samples fails this goal in this challenging situation. For instance, the weighted expected posterior prior in (9) for this situation is

(19) $$\pi_1^{\text{EP}}(\theta) = \sum_{i=1}^{n} \frac{d_i^2}{\|\mathbf{d}\|^2} \frac{d_i}{\sqrt{2\pi}} \exp\left(-\frac{1}{2}(x_i - d_i \theta)^2\right).$$

Although this is, of course, proper, its variance can be shown to be $O(n/\log n)$, so that it becomes increasingly diffuse as $n \to \infty$. Thus the limit is not a stable prior distribution, as one would want.

The problem here is that the training samples corresponding to larger $i$ simply have too little information for them to be useful as training samples. This situation was also encountered in Rodriguez and Pericchi (2001) in dynamic linear models. Their reasonable solution was to only use the most informative training samples to develop intrinsic Bayes factors or expected posterior priors. For instance, a simple modification of (19) would be to truncate the summation at some moderate value $n_0$ (replacing $\|\mathbf{d}\|^2$ by the truncated sum), effectively assigning a "weight" of zero to the low-information training samples. This is an effective option in such situations.

For the general linear model the above phenomenon can also be observed. For clarity we switch to a more standard notation for model selection in the linear model. Suppose for $j = 1, \ldots, q$ that model $M_j$ for the data $\mathbf{Y}$ $(n \times 1)$ is the linear model

$$M_j : \mathbf{Y} = \mathbf{X}_j \boldsymbol{\beta}_j + \boldsymbol{\varepsilon}_j, \qquad \boldsymbol{\varepsilon}_j \sim \mathcal{N}_n(\mathbf{0}, \sigma_j^2 \mathbf{I}_n),$$

where $\sigma_j^2$ and $\boldsymbol{\beta}_j = (\beta_{j1}, \beta_{j2}, \ldots, \beta_{jk_j})'$ are unknown and $\mathbf{X}_j$ is an $(n \times k_j)$ given design matrix of rank $k_j < n$. Let $R_j = |(\mathbf{I} - \mathbf{X}_j(\mathbf{X}_j'\mathbf{X}_j)^{-1}\mathbf{X}_j')\mathbf{y}|^2$ denote the residual sum of squares for $M_j$.

As usual, we utilize the reference prior $\pi_j^N(\boldsymbol{\beta}_j, \sigma_j) = \sigma_j^{-1}$ as the initial noninformative prior. A minimal training sample $\mathbf{y}(l)$, with corresponding design matrix $\mathbf{X}_j(l)$ under $M_j$, is a sample of size $\max\{k_j\} + 1$ such that all $(\mathbf{X}_j'(l)\mathbf{X}_j(l))$ are nonsingular; let $L$ denote the number of such training samples. If $k_j > k_i$,

$$C = \Gamma((n-k_j)/2)\Gamma((k_j - k_i + 1)/2)/[\Gamma((n-k_i)/2)\Gamma(1/2)]$$

and

$$R_j(l) = |(\mathbf{I} - \mathbf{X}_j(l)(\mathbf{X}_j'(l)\mathbf{X}_j(l))^{-1}\mathbf{X}_j'(l))\mathbf{y}(l)|^2,$$

it is shown in Berger and Pericchi ([1996b](#)) that

$$(20) \quad B_{ji}^{\mathrm{A}} = \frac{|\mathbf{X}_i'\mathbf{X}_i|^{1/2}}{|\mathbf{X}_j'\mathbf{X}_j|^{1/2}} \frac{R_i^{(n-k_i)/2}}{R_j^{(n-k_j)/2}} \frac{C}{L} \sum_{l=1}^{L} \frac{|\mathbf{X}_j'(l)\mathbf{X}_j(l)|^{1/2}}{|\mathbf{X}_i'(l)\mathbf{X}_i(l)|^{1/2}} \frac{(R_j(l))^{1/2}}{(R_i(l))^{(k_j - k_i + 1)/2}}.$$

Problems can again arise here if too many of the $|\mathbf{X}_j(l)'\mathbf{X}_j(l)|$ (which are proportional to the "information" in the training samples) are small.

EXAMPLE 16. Consider the special case of testing whether the slope of a linear regression is zero. Thus, let $M_1$ be the model with only the constant term $\beta_1$ and $\mathbf{X}_1' = (1, \ldots, 1)$, and $M_2$ be the model with $(\beta_1, \beta_2)$ and

$$\mathbf{X}_2' = \begin{pmatrix} 1 & \cdots & 1 & 1 & \cdots & 1 & 1 \\ 0 & \cdots & 0 & \delta & \cdots & \delta & 1 \end{pmatrix},$$

with $m = (n-1)/2$ being the number of zeroes and also the number of $\delta$'s. Let $\delta$ be very close to zero. Minimal training samples are then of two types. The high-information minimal training samples are triples $\{y_i, y_j, y_n\}$, where $i \neq j$ range from 1 to $n-1$. There are $m(2m-1)$ such training samples, and they have $|\mathbf{X}_2'(l)\mathbf{X}_2(l)| \cong 2$. The low-information minimal training samples include either one observation from the first $m$ and two observations from the second $m$, or the reverse. There are $m^2(m-1)$ such training samples and they have $|\mathbf{X}_2'(l)\mathbf{X}_2(l)| = 2\delta^2$. Since $\delta$ is very small, the low-information training samples contribute essentially zero to the expression in (20), so that

(with the high-information training samples labelled as $l = 1, \ldots, m(2m - 1)$),

$$B_{21}^{\mathrm{A}} \cong \frac{|\mathbf{X}_1'\mathbf{X}_1|^{1/2}}{|\mathbf{X}_2'\mathbf{X}_2|^{1/2}} \frac{R_1^{(n-1)/2}}{R_2^{(n-2)/2}} \frac{C}{m(m^2+m-1)} \sum_{l=1}^{m(2m-1)} \frac{\sqrt{2}}{\sqrt{3}} \frac{(R_2(l))^{1/2}}{R_1(l)}.$$

As $m$ grows the term involving the training samples clearly goes to zero (since the residual sums of squares for the training samples can be shown to go to nonzero constants as $\delta \to 0$), an undesirable result. Giving equal weight to the (many more) low-information training samples has effectively washed out the effect of the high-information training samples.

The natural solution to this difficulty in the linear model is to weight the training samples according to their information content, that is, choose $p(l) \propto |\mathbf{X}_j(l)'\mathbf{X}_j(l)|$. The problem discussed above will then disappear. Indeed, since there are plenty of high-information training samples available (if $m$ is large), the weighted IBF will have a (nice) intrinsic prior. (This is in contrast to Example 15, where there were not enough high-information training samples to achieve this.) So here weighting works ideally.

The Binet–Cauchy theorem yields the interesting result that

$$p(l) = \frac{|\mathbf{X}_j(l)'\mathbf{X}_j(l)|}{(n-k_j)|\mathbf{X}_j'\mathbf{X}_j|}$$

(i.e., we know the normalization constant for the information-based weighting probabilities), and the weighted IBF then becomes

$$B_{ji}^{\mathrm{A}} = \frac{|\mathbf{X}_i'\mathbf{X}_i|^{1/2}}{|\mathbf{X}_j'\mathbf{X}_j|^{3/2}} \frac{R_i^{(n-k_i)/2}}{R_j^{(n-k_j)/2}} \sum_{l=1}^{L} \frac{C|\mathbf{X}_j'(l)\mathbf{X}_j(l)|^{3/2}}{(n-k_j)|\mathbf{X}_i'(l)\mathbf{X}_i(l)|^{1/2}} \frac{(R_j(l))^{1/2}}{(R_i(l))^{(k_j-k_i+1)/2}}. \quad (21)$$

We do not yet have much experience with use of this IBF, but our current understanding suggests that this will often be better than the usual arithmetic IBF with MTS in linear models. The use of an approximation to a similarly weighted geometric version of the IBF was suggested in de Vos (1993).

Finally, the same issue can be shown to arise with the expected posterior prior in the linear model, so that utilization of the weighted version

$$(22) \qquad \pi_i^{\mathrm{EP}}(\boldsymbol{\beta}_i, \sigma_i^2) = \sum_{l=1}^{L_{\mathrm{M}}} \frac{|\mathbf{X}_i(l)'\mathbf{X}_i(l)|}{(n-k_i)|\mathbf{X}_i'\mathbf{X}_i|} \pi_i^{\mathrm{N}}(\boldsymbol{\beta}_i, \sigma_i^2|\mathbf{y}(l))$$

should be considered.

While the purpose of this paper is not comparison of objective model selection procedures, it is worthwhile to pause and note that the examples we have been considering are challenging for essentially any procedure. As an illustration, consider the most common objective prior used for Bayesian model selection with linear models, the *g-prior*, given by $\pi_i(\sigma_i^2) = 1/\sigma_i^2$ and

$$\pi_i(\boldsymbol{\beta}_i|\sigma_i^2) \qquad \text{is } \mathcal{N}_{k_i}(\mathbf{0}, g\sigma_i^2(\mathbf{X}_i'\mathbf{X}_i)^{-1}).$$

These were proposed in Zellner (1986) for estimation problems. The typical choice of $g$ is $g = n$. Zellner and Siow (1980) suggested a more appropriate (for testing) multivariate Cauchy form for the prior, but it shares with the $g$-prior the underlying scale matrix $\boldsymbol{\Sigma} = n\sigma_i^2(\mathbf{X}_i'\mathbf{X}_i)^{-1}$ which turns out to be quite problematical if it is highly unbalanced.

EXAMPLE 17 (Example 16 continued). Noting that the sample size here is $n$, computation shows that

$$\boldsymbol{\Sigma}_2 = n\sigma_2^2(\mathbf{X}_2'\mathbf{X}_2)^{-1} \cong \sigma_2^2 \begin{pmatrix} 1 & -1 \\ -1 & n \end{pmatrix},$$

so that the information available about $\beta_1$ is vastly different from the information available about $\beta_2$. Indeed, using the $g$-priors with $g = n$ for both $M_1$ and $M_2$ results in the Bayes factor

$$B_{10} = \frac{1}{\sqrt{n+1}} \frac{(\mathbf{y}'\mathbf{y} - n/(n+1)\mathbf{y}'\mathbf{X}_2(\mathbf{X}_2'\mathbf{X}_2)^{-1}\mathbf{X}_2'\mathbf{y})^{-n/2}}{(\mathbf{y}'\mathbf{y} - n/(n+1)\mathbf{y}'\mathbf{X}_1(\mathbf{X}_1'\mathbf{X}_1)^{-1}\mathbf{X}_1'\mathbf{y})^{-n/2}}.$$

For large $n$ and very small $\delta$ [namely, $\delta = \mathrm{o}(n^{-1})$], computation shows that

$$B_{10} \cong \frac{1}{\sqrt{n}} \exp\left(-\frac{(\bar{y} - y_n)^2}{2S^2/n}\right),$$

where $\bar{y}$ and $S^2$ are the usual sample mean and sum of squared deviations. Since the exponential term is bounded in $n$, it follows that $B_{10} \to 0$ as $n$ grows. Hence this Bayes factor is inconsistent under $M_1$, a particularly troubling result.

The difficulty here is that, in a sense, one would like to choose $g = n$ for the information component due to $\beta_1$, but $g = 1$ for the component due to $\beta_2$. The arithmetic IBF and empirical posterior prior (either the weighted or unweighted versions) do this type of adjustment automatically. [It should be mentioned that this would also cause a difficulty with fractional Bayes factors, unless differing fractions are allowed; see De Santis and Spezzaferri (1998a, 1999) and Berger and Pericchi (2001) for discussion.]

In Example 5 it was noted that a problem can also arise with *too informative* training samples, and that it can be wise to restrict attention to training samples whose information content remains modest compared to the information in the entire sample.

EXAMPLE 18 (Example 15 continued). Consider the regression example, but with covariates $d_i = i$. Then the information is rapidly growing with $i$. The expected posterior prior in (19) can then be shown to have variance that is $O(n^{-1})$, so that the prior becomes increasingly (and arguably inappropriately) concentrated as $n \to \infty$. [The same is true if equal weighting is used for the training samples; hence the use of weights in (19) neither helped nor hurt.] Here, simply using only the first, say, $n_0$ training samples (i.e., those with a modest amount of information) would avoid the problem.

It is interesting to note that the common $g$-prior in this situation has variance $n(\sum d_i^2)^{-1} = O(n^{-2})$, which inappropriately concentrates much faster than does the expected posterior prior in (19).

**6. Conclusions.** It is notoriously difficult to develop model selection methodologies that are successful over a wide range of problems. In judging success, our "goal" of developing objective procedures that behave like some reasonable Bayesian procedures may seem to be a rather modest criterion, but it is far stronger than any other criterion we know. We also feel that "testing" a procedure on extreme examples is by far the best method of judging the limits of the procedure, and in suggesting needed refinements. As we have tested intrinsic Bayes factors and expected posterior priors in the years since their development, it has become increasingly clear that the original suggestion—to always use minimal training samples—was too limited. This paper presented a summary of the highlights of these investigations and our suggestions for the needed refinements. The two major conclusions that emerged are:

- In situations, such as censoring, in which certain observations would never be part of an MTS, instead utilize SMTS, which will allow possible involvement of all observations.
- In situations, such as the linear model, in which MTS can contain drastically different information content, consider weighting the training samples (or randomly choosing them) according to their information content.

Random training samples are also useful in other situations, such as when only sufficient statistics, not the actual data, are available. And there are further interesting possibilities that we have not explored, such as forming random training samples by sampling from the data until one has obtained a training sample with at least some pre-specified information content.

Attention in this paper was primarily confined to the arithmetic IBF and the expected posterior prior. However, the generalizations of training samples can (and should) also be used with other training-sample approaches. For instance, the geometric IBF can use the generalizations in exactly the same way as the arithmetic IBF. The median IBF is often preferable to the arithmetic or geometric IBFs from a robustness perspective [Berger and Pericchi (1998)], and randomized training samples can again be utilized directly in its computation. It is not immediately obvious how to utilize weighted training samples with the median IBF, however. The easiest approach is to draw random training samples with probabilities proportional to the weights and then use the median IBF with these training samples.

Finally, it should be noted that this was not meant to be a survey paper, and so we have not dealt with all issues involved in suitably defining training samples. For instance, in Sivaganesan and Lingham (1999) it is shown how transformations of the data are sometimes needed to obtain suitable training samples.

## APPENDIX

LEMMA 2. *In the situation of Example 10, use of the arithmetic IBF based on the Jeffreys-rule prior results in an intrinsic prior with median $O(r^{-1})$.*

PROOF. Since any single observation, censored or uncensored, is an MTS for the Jeffreys-rule prior, (5) leads to the following intrinsic prior:

$$
(23) \qquad \pi^{\mathrm{I}}(\theta) = \pi^{\mathrm{J}}(\theta) \bigg[ \int_0^r \frac{\theta_0 \exp(-\theta_0 x)}{\int \pi^{\mathrm{J}}(\theta) \theta \exp(-\theta x)\, d\theta} \theta \exp(-\theta x)\, dx \\
+ \frac{\exp(-\theta_0 r)}{\int \pi^{\mathrm{J}}(\theta) \exp(-\theta r)\, d\theta} \exp(-r\theta) \bigg].
$$

To study the behavior of the median of this intrinsic prior as $r \to 0$, note that the mass of the first term on the right-hand side of (23) is, switching order of integration,

$$
\iint_0^r \frac{\theta_0 \exp(-\theta_0 x)}{\int \pi^{\mathrm{J}}(\theta) \theta \exp(-\theta x)\, d\theta} \pi^{\mathrm{J}}(\theta) \theta \exp(-\theta x)\, d\theta\, dx \\
= \int_0^r \theta_0 \exp(-\theta_0 x)\, dx \\
= 1 - e^{-\theta_0 r} \to 0 \qquad \text{as } r \to 0.
$$

Hence the median as $r \to 0$ depends only on the second term on the right-hand side of (23). Computation shows that $\int \pi^{\mathrm{J}}(\theta) \exp(-\theta r)\, d\theta \cong 1.5814$ (not depending on $r$). Also, $\exp(-\theta_0 r) \to 1$ as $r \to 0$, so that the median, $\mathrm{med}(r)$, as $r \to 0$ is approximately given by the solution to

$$
0.5 \cong \int_0^{\mathrm{med}(r)} \frac{1}{1.5814} \pi^{\mathrm{J}}(\theta) \exp(-r\theta)\, d\theta.
$$

The change of variables $y = r\theta$, results in the equation

$$0.7907 \cong \int_0^{r \operatorname{med}(r)} y^{-1}(1 - \exp(-y))^{-1/2} \exp(-y)\, dy.$$

Solving this equation for $r \operatorname{med}(r)$ results in the conclusion that $\operatorname{med}(r) \cong 0.191/r$, completing the proof. $\square$

Institute of Statistics and Decision Sciences
Duke University
Durham, North Carolina 27708-0251
USA
e-mail: berger@stat.duke.edu

Departamento de Matemáticas
Universidad de Puerto Rico
P.O. Box 23355
San Juan, Puerto Rico 00931-3355
USA
e-mail: pericchi@cnnet.upr.edu